\documentclass[14pt]{article}
\usepackage{amsmath,amsfonts,amssymb,theorem,verbatim,graphics,epsfig}

\newtheorem{theorem}{Theorem}

 \newtheorem{cor}[theorem]{Corollary}
\def\BEN{\begin{enumerate}}  \def\BI{\begin{itemize}}
\def\EEN{\end{enumerate}}   \def\EI{\end{itemize}}
   \def\sec{\section} 

\def\beq{\begin{eqnarray}} \def\eeq{\end{eqnarray}}

\def\al*#1{\begin{align*}#1\end{align*}}

\def\ga*#1{\begin{gather*}#1\end{gather*}}

\def\alat*#1#2{\begin{alignat*}{#1}#2\end{alignat*}}
\def\bea{\begin{eqnarray*}}
\def\eea{\end{eqnarray*}}
\def\ml*#1{\begin{multline*}#1\end{multline*}}

\def\P{{\mathbb P}}   \def\i{\infty}
\def\E{{\mathbb E}}   \def\R{{\mathbb R}}

\def\te#1{\mathrm{e}^{#1}}

 \def\n{\nu}

\def\w{\omega}  
   
  \def\td{\text{\rm d}}
 
\newtheorem{Thm}{Theorem}

{\theorembodyfont{\normalfont}
 \newtheorem{Def}{Definition}
\newtheorem{Rem}{Remark} 
 
}

\newcommand{\proof}{{\it Proof\ }}

\newcommand{\exit}{{\mbox{\, \vspace{3mm}}}
\hfill\mbox{$\square$}}

\begin{document}

\title{Ruin probability with Parisian delay for a spectrally negative L\'evy risk process}

\author{Irmina Czarna%
\footnote{Department of Mathematics, University of Wroc\l
aw, pl. Grunwaldzki 2/4, 50-384 Wroc\l aw, Poland, e-mail:
czarna@math.uni.wroc.pl}\ \hspace{0.3cm} Zbigniew
Palmowski\footnote{Department of Mathematics, University of Wroc\l
aw, pl. Grunwaldzki 2/4, 50-384 Wroc\l aw, Poland, e-mail:
zbigniew.palmowski@gmail.com}}
\maketitle

\begin{center}
\begin{quote}
\begin{small}
{\bf Abstract.} In this paper we analyze so-called Parisian ruin probability that happens when surplus process stays below zero longer than fixed
amount of time $\zeta>0$. We focus on general spectrally negative L\'{e}vy insurance risk process.
For this class of processes we identify expression for ruin probability in terms of some other quantities
that could be possibly calculated explicitly in many models. We find
its Cram\'{e}r-type and convolution-equivalent
asymptotics when reserves tends to infinity. Finally, we analyze few explicit examples.
\end{small}
\end{quote}
\end{center}

\sec{Introduction}\label{sec:intro}

In risk theory we usually consider classical Cram\'er-Lundberg risk process:
\begin{eqnarray}\label{CL}
 X_t = x + p t -S_t,
\end{eqnarray}
where $x>0$ denotes the initial reserve,
\begin{equation}\label{zlpoisson}
S_t=\sum_{i=1}^{N_t} U_i\end{equation}
is a compound Poisson process. We assume that  $U_i,(i=1,2,...)$ are i.i.d distributed claims (with the distribution function $F$ and tail $\overline{F}:=1-F$ ). Arrival process is a homogeneous Poisson process $N_t$ with intensity $\lambda$.
Premium income is modeled by a constant premium density $p$ and the net profit condition is then $\lambda\nu/p<1$,
where $\E(U_1)=\nu<\infty$.
Lately there has been considered more general setting of a spectrally negative L\'evy process. That is,
$X=\{X_t\}_{t\geq 0}$ is a process with stationary and independent increments with only negative jumps.
We will assume that process starts from $X_0=x$ and
later we will use convention $\P(\cdot|X_0=x)=\P_x(\cdot)$ and $\P_0=\P$.
Such process takes into account not only large claims compensated by a steady income at rate $p>0$, but also
small perturbations coming from gaussian component and possibly additionally (when $\Pi_X(-\infty,0)=\infty$ for jump measure $\Pi_X$ of $X$)
compensated countable infinite number of small claims arriving over each finite time horizon (see e.g. \cite{KyprKlup2}).

One of the most important characteristics in risk theory is a ruin probability defined by
$\P(\tau^{-}_0<\infty)$  for $\tau^{-}_0=\inf\{t>0: X_t < 0 \}$. In this paper we extend this notion to called Parisian ruin probability,
that occurs if process $X$ stays below zero for longer period than a fixed $\zeta>0$.
Formally, we define the excursion below zero:
$$\varsigma_t=\sup\{s<t: \mathbf{1}_{(X_s \geq 0)} \mathbf{1}_{(X_t<0)}=1\}.$$
Parisian time of ruin is given by
$$\tau^{\zeta}=\inf\{t>0: t-\varsigma_t \geq \zeta \}$$
and Parisian ruin probability we define as:
$$\P(\tau^\zeta<\infty|X_0=x)=\P_x(\tau^\zeta<\infty).$$
The case $\zeta=0$ corresponds to classical ruin problem and we do not deal with this case here.
The name for this problem comes from Parisian option that prices
are activated or canceled depending on type of option if
underlying asset stays above or below barrier long enough in a row.
We believe that Parisian ruin probability could be better measure of risk in many situations
giving possibility for insurance company to get solvency. Parisian ruin probability was already considered by Dassios and Wu \cite{DassWu1}
where Parisian ruin probability was found for classical risk process (\ref{CL}) with exponential claims
and for Brownian motion with drift. In Dassios and Wu \cite{DassWu2} Cram\'{e}r-type asymptotics was found for (\ref{CL}).
In this paper using the fluctuation theory we show that these results could be extended
to the case of a general spectrally negative
L\'{e}vy process. In particular we show how to deal with risk process being sum of independent classical risk process (\ref{CL}) and Brownian or $\alpha$-stable motion perturbations.
Additionally, we derive asymptotics of Parisian ruin probability in convolution-equivalent case, that is when
claim size has a heavy-tailed distribution.

This paper is organized as follows. In Section \ref{sec:prel} we introduce basic notions and notations. In Section
\ref{sec:main} we give main representation of Parisian ruin probability.  In Sections \ref{sec:cram} and \ref{sec:conv}
we give asymptotics of Parisian ruin probability in Cram\'{e}r and convolutions-equivalent case, respectively.
Finally, in Section \ref{sec:examples} we analyze
few particular examples.

\sec{Preliminaries}\label{sec:prel}
In this paper we consider a spectrally negative L\'evy process $X=\{X_t\}_{t\geq 0}$, that is a L\'evy process with the L\'evy measure $\Pi_X$
satisfying $\Pi_X(0,\infty)=0$. Assume that $X_0=x>0$ and that $X_t\to\infty$ as $t\to \infty$ a.s.
(that is, reserves of the insurance company increases to infinity a.s.).
This assumption excludes also case of a compound Poisson process with negative jumps.
We assume also that $\E X_1<\infty$. Note that $\E X_1>0$ by drift assumption.
With $X$ we associate the Laplace exponent $\varphi(\beta):=\frac{1}{t}\log \E(e^{\beta X_t})$
defined for all $\beta \geq 0$ and function $\Phi(q)=\sup\{\theta \geq 0: \varphi(\theta)=q\}$.
We will consider also dual process $\widehat{X}_t=-X_t$ which is a spectrally positive L\'evy process with jump measure $\Pi_{\widehat{X}}\left(0,y\right)=\Pi_{X}\left(-y,0\right)$.
Characteristics of $\widehat{X}$ will be indicated by
using a hat over the existing notation for characteristics of $X$.
For process $X$ we define ascending ladder height process $(L^{-1},H)=\{(L^{-1}_t,H_t)\}_{t\geq 0}$:
$$L^{-1}_t:= \left \{ \begin{array}{rl}
\inf\{s>0: L_s>t \} & \textrm{if $t < L_{\infty}$} \\ \infty & \textrm{otherwise}
\end{array} \right.$$
and
$$H_t:= \left \{ \begin{array}{rl}
X_{L_t^{-1}} & \textrm{if $t < L_{\infty}$} \\ \infty & \textrm{otherwise},
\end{array} \right.$$
where $L=\{L_t\}_{t\geq 0}$ is a local time at the maximum (see \cite[p. 140]{Kbook}).
Recall that $(L^{-1},H)$ is a bivariate subordinator with the Laplace exponent
$\kappa(\theta,\beta)=-\frac{1}{t}\log \E\left (e^{-\theta L_t^{-1}-\beta H_t}\textbf{1}_{\{t\leq L_{\infty}\}}\right)$ and with the jump measure
$\Pi_H$.
We define descending ladder height process $(\widehat{L}^{-1},\widehat{H})=\{(\widehat{L}^{-1}_t,\widehat{H}_t)\}_{t\geq 0}$
with the Laplace exponent $\widehat{\kappa}(\theta, \beta)$ constructed from dual process $\widehat{X}$.
Recall that $\widehat{L}_\infty$ has exponential distribution with parameter $\widehat{\kappa}(0,0)$.
Moreover, from the Wiener-Hopf factorization we have:
\begin{equation}\label{kappy}
\kappa(\theta, \beta)=\Phi(\theta)+\beta,\qquad \widehat{\kappa}(\theta, \beta)=\frac{\theta-\varphi(\beta)}{\Phi(\theta)-\beta};
\end{equation}
see \cite[p. 169-170]{Kbook}. Hence
\begin{equation}\label{kappazero}
\widehat{\kappa}(0,0)=\varphi^\prime(0+).
\end{equation}

We introduce a potential measure $\mathcal{U}$ defined by
$$ \mathcal{U}(dx,ds)=\int_0^{\infty} \P(L^{-1}_t \in ds, H_t \in dx) dt$$
with the Laplace transform
$ \int_{[0,\infty)^2} e^{-\theta s -\beta x} \mathcal{U}(dx,ds)=1/\kappa(\theta,\beta) $ and renewal function:
$$U(dx)=\int_{[0,\infty)}\mathcal{U}(dx,ds)=\E\left(\int_0^{\infty} \textbf{1}_{\{H_t\in dx\}} dt \right).$$
For spectrally negative L\'{e}vy process upward ladder height process is
a linear drift and hence the renewal measure is just Lebesgue measure:
\begin{equation}\label{renfunleb}
U(dx)=dx.
\end{equation}
Moreover, from (\ref{kappy}),
\begin{equation}\label{LTdlahatU}
\int_0^\infty e^{-\theta z}\widehat{U}(dz)=\frac{\theta}{\varphi(\theta)};
\end{equation}
see \cite[p. 195]{Kbook}.

We will also use the first passage times:
$$\tau_x^-=\inf\{t\geq 0: X_t<x\},\qquad \tau_x^+=\inf\{t\geq 0: X_t>x\}.$$
Suppose now probabilities $\{\P_{x}\}_{x \in \mathbb{R}}$ corresponds to the
conditional version of $\P$ where $X_{0}=x$ is given. We simply write $%
\P_{0}=\P $. We define also
Girsanov-type change of measure via:
\begin{equation}\label{Girsanov}
\left. \frac{d\P_{x}^{c}}{d\P_{x}}\right| _{\mathcal{F}_{t}}=\frac{\mathcal{E}%
_{t}\left( c\right) }{\mathcal{E}_{0}\left( c\right) }
\end{equation}
for any $c$ for which $\E e^{cX_1}<\infty$, where $\mathcal{E}_{t}\left( c\right) =\exp \{cX_{t}-\psi
\left( c\right) t\}$ is exponential martingale under $\P_x$ and $\mathcal{F}_{t}$ is a natural filtration of $X$.
It is easy to check that under this change of
measure, $X$ remains within the class of spectrally negative processes and
the Laplace exponent of $X$ under $P_{x}^{c}$ is given by
\begin{equation}\label{varphinowe}
\varphi _{c}\left( \theta \right) =\varphi \left( \theta +c\right) -\varphi \left(
c\right)
\end{equation}
for $\theta \geq -c.$
Similarly, by (\ref{kappy}) on $\P^c$ exponent of descending ladder process equals,
\begin{equation}\label{kappynowe}
\widehat{\kappa}_c(0,\beta)=\frac{\varphi \left( \beta +c\right) -\varphi \left(
c\right)}{\beta}
\end{equation}
since on an account of positive drift of $X$ we have $\Phi (0)=0$.

For $q\geq 0$, there exists a function $W^{(q)}: [0,\i) \to [0,\i)$,
called {\it $q$-scale function}, that is continuous and
increasing with the Laplace transform \begin{equation}\label{eq:defW} \int_0^\infty
\te{-\theta x} W^{(q)} (y)  \td y = (\varphi(\theta) - q)^{-1},\quad\theta >
\Phi(q). \end{equation}
We denote $W^{(0)}(x)=W(x)$.
Domain of $W^{(q)}$ is extended to the entire real
axis by setting $W^{(q)}(y)=0$ for $y<0$. For each $y\geq 0$, function $q\to W^{(q)}(y)$ may be analytically
extended to $q\in \mathcal{C}$.
Moreover, let
$$Z^{(q)}(y)=1+q\int_0^yW^{(q)}(z)\,dz.$$
It is known that:
\begin{eqnarray}\label{exit}
\E_x\left(e^{-\theta\tau_y^+}, \tau_y^+<\infty\right)&=&e^{-\Phi(\theta)(y-x)},\\
E_{x}\left( e^{-q\tau _{0}^{-}}\mathbf{1}_{\left( \tau
_{0}^{-}<\infty \right) }\right) &=&Z^{(q)}(x)-\frac{q}{\Phi \left( q\right) }%
W^{(q)}(x)\;,  \label{one-sided-down}
\\
\E_x\left(e^{vX_{\tau^-_0}},\tau_0^- <\infty\right)
&=&e^{vx}\left(Z_v^{(-\varphi(v))}(x) +\frac{\varphi(v)}{\Phi(-\varphi(v))}W_v^{(-\varphi(v))}(x)\right),\label{exit2}
\end{eqnarray}
where $W^{(-\varphi(v))}_v$ and $Z^{(-\varphi(v))}_v$ are scale functions calculated with respect to the measure $\P^v$.
We understand $\theta/\Phi \left( \theta\right) $ in the limiting sense for $%
\theta=0, $ so that
\begin{equation}\label{exit3}
\P_x(\tau _{0}^{-}<\infty )=\left\{
\begin{array}{ll}
1-\varphi ^{\prime }(0)W(x) & \text{if }\varphi ^{\prime }(0)>0 \\
1 & \text{if }\varphi ^{\prime }(0)\leq 0.
\end{array}
\right.
\end{equation}
For details see \cite{Kbook} or \cite[Remark 3]{KypPalm}.
If $\P_x(X_{\tau_0^-}=0, \tau_0^-<\infty)>0$ then we say about downward creeping.
This is possible only when downward ladder process is of unbounded variation with strictly positive drift $d$.
In this case the renewal function $\widehat{U}$ has a strictly positive and continuous density $\widehat{u}$ on $(0,\infty)$
satisfying:
\begin{equation}\label{creeping}
\P_x(X_{\tau_0^-}=0, \tau_0^-<\infty)
=d\widehat{u}(x),
\end{equation}
where $\widehat{u}(0+)=\frac{1}{d}$; see \cite[Theorem 5.9 and Problem 5.5]{Kbook}.

\sec{Main representation}\label{sec:main}
The main representation is given in the next theorem.
\begin{Thm}\label{ThmMain}
Parisian ruin probability for a spectrally negative L\'evy risk process equals:
\begin{eqnarray} \label{PRE}
\lefteqn{\P_x(\tau^\zeta<\infty)=
\P_x(\tau_0^-<\infty)\P(\tau^\zeta<\infty)}\\\nonumber&&+
\left(1-\P(\tau^\zeta<\infty)\right) \int_0^\infty \P(\tau_z^+>\zeta)
\P_x(\tau^{-}_0<\infty, -X_{\tau^{-}_0}\in dz)
\end{eqnarray}
and
\begin{equation}\label{ruinprobref}
\P_x(\tau_0^-<\infty)=\widehat{\kappa}(0,0)\widehat{U}(x,\infty)=1-\varphi^\prime(0+)W(x),
\end{equation}
\begin{eqnarray}\label{ltxuj}
\lefteqn{\int_0^{\infty}e^{-\theta s}\,ds\int_0^\infty \P(\tau_z^+>s)
\P_x(\tau^{-}_0<\infty, -X_{\tau^{-}_0}\in dz)}\\&& =
\frac{\widehat{\kappa}(0,0)\widehat{U}(x,\infty)}{\theta}-\frac{\widehat{\kappa}(0,0)}{\Phi(\theta)}e^{\Phi(\theta)x}
\int_x^\infty e^{-\Phi(\theta)y}\,\widehat{U}(dy)\label{ltx0}
\\ \label{ltx}
&& =\frac{1-\varphi^\prime(0+)W(x)}{\theta}-\frac{1}{\theta}e^{\Phi(\theta)x}\left(Z^{(-\theta)}_{\Phi(\theta)}(x)+\frac{\theta}{\Phi(-\theta)}W^{(-\theta)}_{\Phi(\theta)}(x)\right).
\end{eqnarray}
\end{Thm}
\proof The equation (\ref{PRE}) follows from the strong Markov property and from fact that $X$ has only negative jumps.
Moreover, identity
\begin{eqnarray}\label{Dualequi}
\P_x(\tau^-_0<\infty)=\P(\widehat{\tau}^+_x<\infty)
\end{eqnarray}
and observation that for $x>0$:
\begin{eqnarray}\label{tau0}
\P(\widehat{\tau}^+_x<\infty)=\widehat{\kappa}(0,0)\widehat{U}(x,\infty)
\end{eqnarray}
complete proof of (\ref{ruinprobref}) in view of (\ref{exit3}) (for (\ref{tau0}) see \cite[p. 187]{Kbook}).
Finally, note that from (\ref{exit}):
\begin{equation}\label{LapXeq}\int_0^{\infty}e^{-\theta s}\P(\tau_z^+>s)\,ds = \frac{1}{\theta}\left(1- \E\left(e^{-\theta\tau_z^+}\mathbf{1}_{\left(\tau_z^+<\infty\right) }\right)\right)
=\frac{1}{\theta}\left(1- e^{-\Phi(\theta)z}\right).
\end{equation}
Further, in case of downward creeping, distribution
$\P_x(\tau^{-}_0<\infty, -X_{\tau^{-}_0}\in dz)$ might have an atom at $0$ but it does not have influence on the rhs of
(\ref{ltx}):
\begin{eqnarray}\lefteqn{\int_0^{\infty}e^{-\theta s}\,ds\int_0^\infty \P(\tau_z^+>s)
\P_x(\tau^{-}_0<\infty, -X_{\tau^{-}_0}\in dz)}\nonumber\\&&=
\frac{1}{\theta}\int_0^\infty \left(1- e^{-\Phi(\theta)z}\right)\P_x(\tau^{-}_0<\infty, -X_{\tau^{-}_0}\in dz)
\nonumber\\&&=\frac{\P_x(\tau_0^-<\infty)}{\theta}-\frac{1}{\theta}\E_x\left[ e^{\Phi(\theta)X_{\tau_0^-}}, \tau^{-}_0<\infty\right].\label{waznarep}
\end{eqnarray}
Identity (\ref{exit2}) completes proof of (\ref{ltx}). The equation (\ref{ltx0}) follows from following equality:
\begin{eqnarray} \label{Dual}
\P_x(-X_{\tau^{-}_0}\in dz,\tau^{-}_0<\infty)=\P(\widehat{X}_{\widehat{\tau}^+_x}-x \in dz,\widehat{\tau}^{+}_x<\infty )
\end{eqnarray}
and \cite[Problem 5.5]{Kbook}.

\exit

\begin{Rem}
Note that
$$\int_0^\infty e^{-\beta x}W^{(-\theta)}_{\Phi(\theta)}(x)\,dx=\varphi(\beta+\Phi(\theta))^{-1}, \quad Z^{(-\theta)}_{\Phi(\theta)}(x)=1-\theta\int_0^xW^{(-\theta)}_{\Phi(\theta)}(y)\,dy.$$
\end{Rem}

\begin{Rem}\label{Brown}
Consider now particular case of spectrally negative L\'evy process:
\begin{equation}\label{reprfinalowa}
X_t=x+pt-S_t+\sigma B_t,\end{equation}
where $p>0$, $\sigma\geq 0$, $S_t$ is a subordinator of bounded variation with jump measure $\Pi_{\widehat{X}}$ and $B_t$ is a Brownian motion independent of $S_t$.
Moreover, since $S_t$ has bounded variation, then:
$$\nu_0=-\int_{-\infty}^0z\Pi_X(dz)=\int_0^{\infty}\overline{\Pi}_{\widehat{X}}(z)\,dz<\infty,$$
where $\overline{\Pi}_{\widehat{X}}(z)=\Pi_{\widehat{X}}(z,\infty)$.
Note that $\rho_0=\E S_1/p<1$ because $\E X_1>0$.
We recall now Pollaczek-Khintchine formula (see \cite[Theorem 3.1]{Huzak2}):
\begin{equation}\label{pch}
\P_x(\tau_0^-<\infty)=\widehat{\kappa}(0,0)\widehat{U}(x,\infty)=1-(1-\rho_0)\sum_{n=0}^\infty\rho_0^n\left(K^{(n+1)*}*M^{n*}\right)(x),
\end{equation}
where \begin{equation}\label{pihatH}M(dz)=\frac{1}{\nu_0}\overline{\Pi}_{\widehat{X}}(z)\,dz,\end{equation}
and when $\sigma >0$:
$$\int_0^\infty e^{-\theta z} K(dz)=\frac{p\theta}{p\theta+\frac{1}{2}\theta^2\sigma^2},$$
hence $K(dz)=\frac{2p}{\sigma^2}e^{-(2pz)/\sigma^2}\,dz$.

{\bf The case of $\sigma=0$.} If there is no Gaussian component, that is $X_t$ is drift process minus subordinator of bounded variation,
then $K(dz)=\delta_0(dz)$.
The equation (\ref{pch}) allows also to identify $\int_0^{\infty}e^{-\theta s}\,ds\int_0^\infty \P(\tau_z^+>s)
\P_x(\tau^{-}_0<\infty, -X_{\tau^{-}_0}\in dz)$ using (\ref{pch}), (\ref{ruinprobref}) and (\ref{ltx0}).
Moreover, 
\begin{equation}\label{ro}
\rho=\P(\tau_0^-<\infty)=1-\frac{\E X(1)}{p}=\frac{\E S_1}{p}=\rho_0.\end{equation}

\end{Rem}

To identify ruin probability $\P_x(\tau^\zeta<\infty)$ we need to find constant \linebreak $\P(\tau^\zeta<\infty)$, which
is given in the next theorem.\\

 Denote by $p^+(s)=\P_\epsilon(\tau_0^-<s)$ probability that the excursion above $0$ is shorter than $s$ and let
\begin{eqnarray*}
p(s,t)
&=&\int_0^\infty \P(\tau^+_{z+\epsilon}\leq t)\,\P_\epsilon(\tau_0^-<s, -X_{\tau_0^-}\in dz)
\\&&
+\P(\tau^+_{\epsilon}\leq t)
\P_\epsilon(\tau_0^-<s, X_{\tau_0^-}=0)
\end{eqnarray*} be probability that upper excursion above $0$ is shorter than $s$ and the first consecutive excursion below $0$,
which is shifted downward by $-\epsilon$,
is shorter then $t$. Note that $p^+(s)=p(s, \infty)$.

\begin{Thm}\label{ThmMain2}
\begin{description}\item{(i)}
If $X$ is a process of bounded variation, then
\begin{eqnarray}\label{P0}
\nonumber \P(\tau^\zeta<\infty)&=&\frac{\int_0^\infty \P(\tau_z^+>\zeta)
\P(\tau^{-}_0<\infty, -X_{\tau^{-}_0}\in dz)}{1- \rho+\int_0^\infty \P(\tau_z^+>\zeta)
\P(\tau^{-}_0<\infty, -X_{\tau^{-}_0}\in dz)},
\end{eqnarray}
where
\begin{eqnarray}\nonumber
\lefteqn{\int_0^{\infty}e^{-\theta s}\,ds\int_0^\infty \P(\tau_z^+>s)
\P(\tau^{-}_0<\infty, -X_{\tau^{-}_0}\in dz)}\\
 &&= \frac{1}{\theta p} \int_0^{\infty}\left(1- e^{-\Phi(\theta)z}\right) \overline{\Pi}_{\widehat{X}}(z) dz.
\end{eqnarray}
\item{(ii)}
If $X$ is a process of unbounded variation, then
\begin{equation}\label{regular}
\P(\tau^\zeta<\infty)=\lim_{b\to\infty}\lim_{\epsilon\downarrow 0}\frac{p^+(b)-p(b,\zeta)}{1-p(b,\zeta)}.
\end{equation}
\end{description}
\end{Thm}
\proof
If $X$ is a process of bounded variation then by \cite[Cor. VII.5]{bertbook} $0$ is irregular for $(-\infty,0)$.
Since by drift assumption we excluded the case of compound Poisson process $X$, it means that $0$ is regular for $(0,\infty)$.
Now point (i) follows from Theorem \ref{ThmMain} by taking $x=0$ and by using
identities (\ref{ro}) and  (\ref{BoundVar}).
In the proof of (ii) we adapt the idea from the proofs of main results of Dassios and Wu \cite{DassWu1} and Dassios and Wu \cite{DassWu3}.\\
Let $\delta_0^\pm=0$. First we will define the sequence of stopping times:
\begin{eqnarray*}
\sigma_n^+=\inf\{t>\delta^+_n: X_t\leq -\epsilon\},&&\qquad
\sigma_n^-=\inf\{t>\delta^-_n: X_t=\epsilon\},\\
\delta_{n+1}^+=\inf\{t>\sigma_n^+: X_t=0\},&&\qquad \delta_{n+1}^-=\inf\{t>\sigma_n^-: X_t \leq 0\}
\end{eqnarray*}
and the processes:
$$X^{\pm}_t=\left\{\begin{array}{lr}
X_t\pm\epsilon&\mbox{if $\delta_n^\pm\leq t<\sigma_n^\pm$}\\
X_t&\mbox{if $\sigma_n^\pm\leq t<\delta_{n+1}^\pm$}.
\end{array}\right.
$$
Moreover, let $b>0$ and
$$\tau_b=\inf\{t>0: t-\sup\{s<t: \mathbf{1}_{(X_s \leq 0)} \mathbf{1}_{(X_t>0)}=1\}>b\}.$$
Similarly, we define stopping times $\tau_b^\pm$ and $\tau^\zeta_\pm$ for $X^\pm$.
Observe that:
\begin{equation}\label{estimation}
\P(\tau^\zeta_+<\tau_b^+)\leq \P(\tau^\zeta<\tau_b)\leq \P(\tau^\zeta_-<\tau_b^-).
\end{equation}
We will find $\P(\tau^\zeta_+<\tau_b^+)$ first by decomposing the path of the process $X^+$ into excursions above and below $0$.
Formally, probability of event $A_j$, that the first excursion below $0$
of length greater than $\zeta$ is $j$th excursion and it happens before the first excursion above $0$ longer than $b$,
equals:
$$\P(A_j)=p(b,\zeta)^{j-1}(p^+(b)-p(b,\zeta)).$$
Summing over $j=1,2,\ldots$ gives:
\begin{equation}\label{plus}
\P(\tau^\zeta_+<\tau_b^+)=\frac{p^+(b)-p(b,\zeta)}{1-p(b,\zeta)}.
\end{equation}
Similarly,
\begin{equation}\label{minus}
\P(\tau^\zeta_-<\tau_b^-)=\P(\tau_\epsilon^+>\zeta)+\P(\tau_\epsilon^+\leq \zeta)\frac{p^+(b)-p(b,\zeta)}{1-p(b,\zeta)}.
\end{equation}
Recall that 
process of bounded variation $X$
was excluded from our considerations.
Thus by regularity of $0$ for $(-\infty,0)$ it follows that $0$ is regular for $(0,\infty)$ (see \cite[Theorem 6.5, p. 142]{Kbook}).
Straightforward consequence of this fact is that $\lim_{\epsilon\downarrow 0}\P(\tau_\epsilon^+\leq \zeta)=1$
and that
$\P(\tau^\zeta_-<\tau_b^-)$ and $\P(\tau^\zeta_+<\tau_b^+)$ have the same limits as $\epsilon\downarrow 0$.
From (\ref{estimation}) we then derive the assertion of the theorem.
\exit

\begin{Rem}
To find $\P(\tau^\zeta<\infty)$ for the process (\ref{reprfinalowa}) with $\sigma=0$ one can use the following identity for $z>0$
(see \cite[Corollary 4.5]{Huzak2} and \cite[p. 105 and Cor. 7.5]{Kbook}):
\begin{eqnarray}\label{BoundVar}
\P(-X_{\tau^{-}_0}\in dz, \tau^{-}_0<\infty)=\rho_0\Pi_{\widehat{H}}(dz)=\frac{1}{p} \overline{\Pi}_{\widehat{X}}(z) dz
\end{eqnarray}
and then
\begin{eqnarray}\nonumber
\lefteqn{\int_0^{\infty}e^{-\theta s}\,ds\int_0^\infty \P(\tau_z^+>s)
\P_x(\tau^{-}_0<\infty, -X_{\tau^{-}_0}\in dz)}\\\nonumber
 &&=\frac{1}{\theta p}\int_0^{\infty}\left(1- e^{-\Phi(\theta)z}\right)\overline{\Pi}_{\widehat{X}}(z) dz.
\end{eqnarray}

\end{Rem}

\begin{Rem}\label{Brown2}
The probability $1-p(s, t)$ could be identified using double Laplace transform.
Indeed, by \cite[Ex. 6.7 (i), p. 176]{Kbook} and (\ref{exit}),
\begin{eqnarray*}\beta\w \int_0^\infty\int_0^\infty (1-p(s,t))e^{-\beta t} e^{-\w s}\,dt\,ds
&=&1-e^{-\Phi(\beta)\epsilon}\E_\epsilon \left( e^{-\Phi(\beta)X_{\tau_0^-} -\w  \tau_0^-}\right) \\&=&
1-\frac{\E \left(e^{\Phi(\beta)\underline{X}_{e_\w}}\mathbf{1}_{(-\underline{X}_{e_\w}>\epsilon)}\right)}{\E e^{\Phi(\beta)\underline{X}_{e_\w}}},
\end{eqnarray*}
where $e_\w$ is independent of $X$ exponential random variable with intensity $\w$.
Moreover, from  \cite{KypPalm} we know that
\begin{equation}\label{inflapl}
\E e^{\Phi(\beta)\underline{X}_{e_\w}}
=\frac{\w(\Phi(\beta)
-\Phi \left( \w\right) )}{\Phi \left( \w\right) (\beta-\w)
}
\end{equation}
and that
\begin{equation}\label{infdensity}P(-\underline{X}_{e_\w}\in dz)=\frac{\w}{\Phi(\w)}W^{(\w)}(dz)-\w W^{(\w)}(z)\,dz.\end{equation}
Hence choosing function $n(\epsilon)$ such that limit
\begin{equation}\label{mw}m(\w)=\lim_{\epsilon\downarrow 0}\frac{P(-\underline{X}_{e_\w}\leq \epsilon)}{n(\epsilon)}\end{equation}
exists and is finite, we have:
\begin{equation}\label{mainrepras}
\lim_{\epsilon\downarrow 0}\frac{1-p(s,t)}{n(\epsilon)}=q(s,t),
\end{equation}
where
\begin{equation}\label{mainrepras2}
\int_0^\infty\int_0^\infty q(s,t)\,dt\,ds=\frac{m(\w)\Phi \left( \w\right) (\beta-\w)}{\beta\w^2(\Phi(\beta)
-\Phi \left( \w\right) )}.
\end{equation}
Then also we have:
\begin{equation}\label{newpzetazero}
\P(\tau^\zeta<\infty)=\lim_{b\to\infty}\frac{q(b,\zeta)-q(b,\infty)}{q(b,\zeta)}.
\end{equation}
\end{Rem}

\section{Cram\'er's estimate}\label{sec:cram}

In this section we derive exponential asymptotics of Parisian ruin probability.
Assume that there exists $\gamma>0$ satisfying $\widehat{\varphi}(\gamma)=\varphi(-\gamma)=0$ and that
$\widehat{\varphi}(\theta)$ is finite in the neighborhood of $\gamma$.
Then $\E e^{-\gamma X_1}<\infty$ and we can define new measure $\P^{-\gamma}$ via
(\ref{Girsanov}).
Define $\widehat{U}_{\gamma}(dx)=\widehat{U}_{\gamma}^{(0)}(dx):=e^{\gamma x} \widehat{U}(dx)$ and
$$\mu = \int_{0}^{\infty} x \widehat{U}_{\gamma}^{(1)}(dx),$$ where
$\widehat{U}_{\gamma}^{(q)}(dx)=\int_{0}^{\infty} e^{-(qt+\gamma x)}\P(\widehat{H}_t\in dx)dt$
for $q \geq 0$.
Note from \cite{BertDon} that $\widehat{U}_{\gamma}$ is a renewal function of the ladder height process calculated on $\P^{-\gamma}$.
Moreover, from \cite{BertDon} we have:
\begin{equation}\label{murazjeszcze}
\mu=\frac{\partial \widehat{\kappa}_{-\gamma}(0,\beta)}{\partial \beta}|_{\beta=0}=\frac{\partial \widehat{\kappa}(0,\beta)}{\partial \beta}|_{\beta=-\gamma}.
\end{equation}
Note also that drift of $\widehat{X}$ on $\P^{-\gamma}$ is positive, since by (\ref{varphinowe}) its Laplace exponent equals then
$$\widehat{\varphi}_{-\gamma}(\theta)=\varphi(-\gamma-\theta)$$
and hence $\widehat{\varphi}_{-\gamma}^\prime(0+)=-\varphi^\prime(-\gamma-)>0$. Thus
\begin{equation}\label{zmdryfu}
\P^{-\gamma}(\widehat{\tau}_x^+<\infty)=1.
\end{equation}

\begin{Thm}\label{Crameras}
We assume Cram\'er conditions, that is, that there exists a $\gamma \in (0,\infty)$ such that $\widehat{\varphi}(\gamma)=0$ and that
$\widehat{\varphi}(\theta)$ is finite in the neighborhood of $\gamma$. We also assume that support of $\Pi$ is not
lattice when $\widehat{\Pi}(R)<\infty$. We have:
\begin{eqnarray} \label{cas}
\lim_{x \uparrow \infty}e^{\gamma x} \P_x(\tau^\zeta<\infty)
 &=& \P(\tau^\zeta<\infty)\frac{\widehat{\kappa}(0,0)}{\gamma\mu}+ (1-\P(\tau^\zeta<\infty))f^{(c)}(\zeta),
\end{eqnarray}
where
$$\int_0^\infty e^{-\theta s}f^{(c)}(s)\;ds=\frac{\widehat{\kappa}(0,0)}{\gamma\mu\theta}-\frac{1}{(\gamma+\Phi(\theta))^2\mu}
$$
and $\P(\tau^\zeta<\infty)$ is given in Theorem \ref{ThmMain2}.
If $\mu=\infty$, then LHS of (\ref{cas}) is understood to be 0.
\end{Thm}

\proof
By (\ref{tau0}) and the Key Renewal theorem saying that
$\widehat{U}_{\gamma}(dx)$  on $(0,\infty)$ converges weakly as a measure to $\mu^{-1} dx$
(see \cite[p. 188]{Kbook} and \cite{BertDon}) we have:
\begin{eqnarray}
\nonumber\lim_{x \uparrow \infty}e^{\gamma x}\P(\widehat{\tau}^+_x<\infty)&=&\widehat{\kappa}(0,0)\lim_{x \uparrow \infty}e^{\gamma x}
\widehat{U}(x,\infty)=\widehat{\kappa}(0,0)\lim_{x \uparrow \infty}\widehat{U}_{\gamma}(x,\infty)
=\frac{\widehat{\kappa}(0,0)}{\gamma\mu}. \label{Krthm}
\end{eqnarray}
Moreover, by (\ref{waznarep}) and (\ref{Dual}),
\begin{eqnarray}\nonumber
\lefteqn{\int_0^{\infty}e^{-\theta s}\,ds\int_0^\infty \P(\tau_z^+>s)
\P_x(\tau^{-}_0<\infty, -X_{\tau^{-}_0}\in dz)}
\\
\nonumber &&=
\frac{\P(\widehat{\tau}^+_x<\infty)}{\theta}-\frac{1}{\theta}\E\left[ e^{-\Phi(\theta)(\widehat{X}_{\widehat{\tau}_x^+}-x)}, \widehat{\tau}^{+}_x<\infty\right]
\end{eqnarray}
and by (\ref{zmdryfu}), the Optional Stopping Theorem and \cite[Problem 5.5]{Kbook} we have:
\begin{eqnarray*}
\lefteqn{\lim_{x \uparrow \infty}e^{\gamma x}\E\left[ e^{-\Phi(\theta)(\widehat{X}_{\widehat{\tau}_x^+}-x)}, \widehat{\tau}^{+}_x<\infty\right] =
\lim_{x \uparrow \infty} \E^{-\gamma}\left[ e^{-(\Phi(\theta)+\gamma)(\widehat{X}_{\widehat{\tau}_x^+}-x)}, \widehat{\tau}^{+}_x<\infty\right]}\\&&
=\lim_{x \uparrow \infty} \E^{-\gamma}\left[ e^{-(\Phi(\theta)+\gamma)(\widehat{X}_{\widehat{\tau}_x^+}-x)}\right]=
\widehat{\kappa}_{-\gamma}(0,\Phi(\theta)+\gamma)\lim_{x \uparrow \infty} \int_x^\infty e^{-(\Phi(\theta)+\gamma)(y-x)}\,\widehat{U}_\gamma(dy)
\\&&=\widehat{\kappa}_{-\gamma}(0,\Phi(\theta)+\gamma)\lim_{x \uparrow \infty} \int_0^\infty e^{-(\Phi(\theta)+\gamma)y}\,\widehat{U}_\gamma(x+dy)=
\frac{\varphi(\Phi(\theta))}{(\gamma+\Phi(\theta))^2\mu}\\&&=\frac{\theta}{(\gamma+\Phi(\theta))^2\mu}.
\end{eqnarray*}
This completes the proof in view of Theorem \ref{ThmMain}.

\exit

\section{Convolution-equivalent case }\label{sec:conv}
We start from introducing class of L\'evy processes, which we will deal with in this section.
Process $\widehat{X}$ considered here
have jump measure $\Pi_{\widehat{X}}$ in a class $\mathcal{S}^{(\alpha)}$ of convolution equivalent distributions,
which could be formally defined as follows:
\begin{Def} (Class $\mathcal{L}^{(\alpha)}$)
For a parameter $\alpha \geq 0$ we say that distribution function G on $[0,\infty)$ with tail $\overline{G}=1-G$
belongs to class $\mathcal{L}^{(\alpha)}$ if
\begin{itemize}
\item[$(i)$]  $\overline{G}(x)>0$ for each $x \geq 0$,
\item[$(ii)$] $\lim_{u \rightarrow \infty} \frac{\overline{G}(u-x)}{\overline{G}(u)}=e^{\alpha x} \textrm{   for each $x \in \R$, and G is nonlattice}$,
\item[$(iii)$] $\lim_{n \rightarrow \infty} \frac{\overline{G}(n-1)}{\overline{G}(n)}=e^{\alpha} \textrm{   if G is lattice}$ (then assumed of span $1$).
\end{itemize}
\end{Def}

\begin{Def} (Class $\mathcal{S}^{(\alpha)}$)
We say that G belongs to class $\mathcal{S}^{(\alpha)}$
if
\begin{itemize}
\item[$(i)$] $G \in \mathcal{L}^{(\alpha)}$ for some $\alpha \geq 0$,
\item[$(ii)$] for some $M_0<\infty$, we have
\begin{eqnarray}
\lim_{u \rightarrow \infty} \frac{\overline{G^{*2}}(u)}{\overline{G}(u)}=2M_0,
\end{eqnarray}
where $\overline{G^{*2}}(u)=1-G^{*2}(u)$ and $*$ denote convolution.
\end{itemize}
\end{Def}
For all $a \in \R$ such that the following integral is finite we define the moment generating function $\delta$ such that
$$\delta_a(G)=\int_0^{\infty} e^{au} G(du),$$
where G is a finite distribution function.

Recall that $\widehat{X}=-X$ is a spectrally positive L\'evy process. Throughout this section we assume that
for $\widehat{X}$ we have:
\begin{itemize}
\item[$(i)$] \begin{eqnarray}\label{Con1}
\overline{\Pi}_{\widehat{X}} \in \mathcal{S}^{(\alpha)} \textrm{ for } \alpha > 0
\end{eqnarray}
and
\begin{eqnarray}\label{Con1b}
\int_0^x\overline{\Pi}_{\widehat{X}}(y)\;dy \in \mathcal{S}^{(0)}
\end{eqnarray}
\item[$(ii)$]
\begin{eqnarray}\label{Con2b} \widehat{\varphi}(\alpha)<0\qquad \textrm{for}\quad \alpha>0;\end{eqnarray}
\item[$(iii)$]
\begin{eqnarray}\label{Con2} e^{-q}\delta_{\alpha}(\widehat{H})<1,\quad \textrm{where}\quad q=\lim_{\beta \downarrow 0} \frac{-\widehat{\varphi}(-\beta)}{\kappa(0,-\beta)}\quad \textrm{and}\quad \delta_{\alpha}(\widehat{H})=\delta_{\alpha}(\widehat{H}_1). \end{eqnarray}
\end{itemize}
The first condition gives:
\begin{eqnarray}\label{Con1c}
\overline{\Pi}_{\widehat{H}} \in \mathcal{S}^{(\alpha)} \textrm{ for } \alpha \geq  0;
\end{eqnarray}
Latter condition (iii) has a force when $\alpha>0$; by drift assumption for $\alpha=0$ this condition is automatically satisfied.

We define relation: $f(x)\sim g(x)$  iff  $\lim_{x \rightarrow \infty} f(x)/g(x)=1$.

\begin{Thm}\label{Thmconveq}
Under assumptions (\ref{Con1})-(\ref{Con2}) asymptotic Parisian ruin probability equals:
\begin{eqnarray} \nonumber
 \P_x(\tau^\zeta<\infty) \sim \E X_1 \left(\frac{\alpha}{\widehat{\varphi}(\alpha)}\right)^2 \left(\P(\tau^\zeta<\infty)+
(1-\P(\tau^\zeta<\infty))f^{(e)}(\zeta)\right)\int_x^\infty\overline{\Pi}_{\widehat{X}}(y)\,dy,
\end{eqnarray}
where
\begin{equation}\label{fe}\int_0^\infty e^{-\theta s}f^{(e)}(s)\,ds=\frac{1}{\theta}
\int_0^{\infty}\left(1- e^{-\Phi(\theta)z}\right) B(z)dz
\end{equation}
and
\begin{eqnarray*}
B(z)
= \frac{e^{-\alpha z}}{E X_1}\left(-\widehat{\varphi}(\alpha)+\alpha \int_z^\infty e^{\alpha y}
\overline{\Pi}_{\widehat{X}}(y)\,dy\right).\end{eqnarray*}
Above for $\alpha=0$ the quantity $-\widehat{\varphi}(\alpha)/\alpha$ is understood in limiting sense and equals $-\widehat{\varphi}^\prime(0+)=EX_1$.
\end{Thm}
\proof
From \cite[Theorem 6.2]{KyprKlup} we have that:
\begin{eqnarray}\label{NCLap}
\P(\widehat{\tau}^+_x<\infty)=\widehat{\kappa}(0,0)\widehat{U}(x,\infty)
\sim EX_1 \left(\frac{\alpha}{\widehat{\varphi}(\alpha)}\right)^2\int_x^\infty\overline{\Pi}_{\widehat{X}}(y)\,dy
\end{eqnarray}
and that:
$$
\lim_{x \uparrow \infty} \P(\widehat{X}_{\widehat{\tau}^+_x}-x \in dz|\widehat{\tau}^{+}_x<\infty)= B(z)dz,
$$
where $$
B(z)=-\frac{d}{dz}\frac{e^{-\alpha z}}{E X_1}\left(-\frac{\widehat{\varphi}(\alpha)}{\alpha} +\int_z^\infty \left(e^{\alpha y}-e^{\alpha z}\right)
\overline{\Pi}_{\widehat{X}}(y)\,dy\right).
$$
The equations (\ref{LapXeq}) and (\ref{Dual}), dominated convergence theorem and Theorem \ref{ThmMain} complete the proof.

\exit

\section{Examples}\label{sec:examples}

\subsection{General classic risk process (\ref{CL})}
For the process (\ref{CL}) we have
$\Pi_{\widehat{X}}(dz)=\lambda F(dz)$ and $\varphi(\theta)=p\theta-\lambda +\lambda\int_0^\infty e^{-\theta z}\,F(dz)$.\\
By (\ref{BoundVar}),
$$ \P(\tau^{-}_0<\infty, -X_{\tau^{-}_0}\in dz)=\frac{\lambda}{p} \overline{F}(z)\,dz,$$
and hence
\begin{eqnarray}\nonumber
\lefteqn{\int_0^{\infty}e^{-\theta s}\,ds\int_0^\infty \P(\tau_z^+>s)
\P(\tau^{-}_0<\infty, -X_{\tau^{-}_0}\in dz)}\\
 &&= \frac{\lambda}{\theta p} \int_0^{\infty}\left(1- e^{-\Phi(\theta)z}\right) \overline{F}(z)\, dz.\label{ltcomp}
\end{eqnarray}
{\bf Probability $\P(\tau^\zeta<\infty)$ could be found} using Theorem \ref{ThmMain2} with (\ref{ltcomp}) and
$$\rho=\P(\tau^{-}_0<\infty)=\frac{\lambda \nu}{p},$$
where $\n=\int_0^\infty y\,F(dy)$.
Theorem \ref{ThmMain} gives then {\bf Parisian ruin probability $\P_x(\tau^\zeta<\infty)$ for all $x\geq 0$}.

To find {\bf Cram\'er asymptotics},
note that
$\Phi(\theta)$ and $\gamma$ solve the equations:
$$\int_0^\infty e^{-\Phi(\theta) z}\,F(dz)=(\lambda-p\Phi(\theta)+\theta)/\lambda,\qquad
\int_0^\infty e^{\gamma z}\,F(dz)=(\lambda+p\gamma)/\lambda.$$
Moreover, from (\ref{kappazero}) we have:
\begin{eqnarray}\label{kapa}
\widehat{\kappa}(0,0)=
\varphi^\prime(0+)=p-\lambda \n.
\end{eqnarray}
By (\ref{murazjeszcze}) using fact that $\varphi(-\gamma)=0$ we have:
$$\mu=\lambda\int_{0}^{\infty} ye^{\gamma y} \Pi_{\widehat{X}}(y,\infty)\,dy
=\lambda\int_{0}^{\infty}y e^{\gamma y} \overline{F}(y)\,dy
$$
and hence
\begin{eqnarray*}
\lim_{x \uparrow \infty}e^{\gamma x} \P_x(\tau^\zeta<\infty)
 &=& \P(\tau^\zeta<\infty)\left (\frac{\lambda \gamma }{p-\lambda \nu}\int_{0}^{\infty} y e^{\gamma y} \overline{F}(y) \,dy \right)^{-1}\\&&+ (1-\P(\tau^\zeta<\infty))f^{(c)}(\zeta),
\end{eqnarray*}
where
\begin{eqnarray*}\int_0^\infty e^{-\theta s}f^{(c)}(s)\;ds&=&
\left(\lambda\int_{0}^{\infty}y e^{\gamma y} \overline{F}(y)\,dy\right)^{-1}
\left(\frac{p-\lambda \n}{\gamma\theta}-
\frac{1}{(\gamma+\Phi(\theta))^2}\right).
\end{eqnarray*}

If $F\in \mathcal{S}^{(\alpha)}$ for $\alpha \geq 0$ then
$\P(\tau^\zeta<\infty)$ is given in Theorem \ref{Thmconveq}.\\

\subsection{Classic risk process (\ref{CL}) with exponential jumps}
\begin{cor}
Assume that $X_t$ is Cram\'er-Lundberg risk process (\ref{CL}) with exponential claims $F(dz)=\xi e^{-\xi z}\,dz$, where $\xi=1/\nu$.
Then
\begin{eqnarray}\label{expo}
 \P_x(\tau^\zeta<\infty)= \frac{\lambda}{p\xi}e^{-(\frac{p\xi-\lambda}{p})x}\left( \frac{p\xi D}{p\xi -\lambda(1-D) }\right)
\end{eqnarray}
where $$D=1-\int_0^\zeta \sqrt{\frac{p\xi}{\lambda}}e^{-(\lambda+p\xi)t}t^{-1}I_1(2t\sqrt{p\lambda\xi})dt$$
and $I_1(x)$ is modified Bessel function of the first kind.
\end{cor}
\proof
From (\ref{ruinprobref}) and (\ref{eq:defW}) (see also \cite[p. 63]{Assbook}) we have:
\begin{equation}\label{exporuin}
\widehat{\kappa}(0,0)\widehat{U}(x,\infty)=\frac{\lambda}{p\xi}e^{-\left(\frac{p\xi-\lambda}{p}\right)x}.
\end{equation}
Moreover, by lack-of-memory property of exponential distribution, distribution of undershoot of $0$ is also
exponential with intensity $\xi$. Hence:
\begin{eqnarray}\label{Lap1}
 \int_0^{\infty}e^{-\theta s}\,ds \int_0^\infty \P(\tau_z^+>s)\P_x(\tau^{-}_0<\infty, -X_{\tau^{-}_0}\in dz)
=\frac{\lambda}{p\xi}e^{-\left(\frac{p\xi-\lambda}{p}\right)x}\frac{\Phi(\theta)}{\theta(\Phi(\theta)+\xi)}
\end{eqnarray}
and by (\ref{ltcomp})
\begin{eqnarray}\label{Lap2}
 \int_0^{\infty}e^{-\theta s}\,ds \int_0^\infty \P(\tau_z^+>s)\P(\tau^{-}_0<\infty, -X_{\tau^{-}_0}\in dz)
=\frac{\lambda}{p\xi}\frac{\Phi(\theta)}{\theta(\Phi(\theta)+\xi)}.
\end{eqnarray}
Inverting the Laplace transforms (\ref{Lap1}) and (\ref{Lap2}) with respect to $\theta$ (see \cite{Bete}) gives:
\begin{eqnarray}
\int_0^\infty \P(\tau_z^+>\zeta)\P_x(\tau^{-}_0<\infty, -X_{\tau^{-}_0}\in dz)=
\frac{\lambda}{p\xi}e^{-(\frac{p\xi-\lambda}{p})x}D\label{exponadskok}
\end{eqnarray}
and
\begin{eqnarray}
\int_0^\infty \P(\tau_z^+>\zeta)\P(\tau^{-}_0<\infty, -X_{\tau^{-}_0}\in dz)=
\frac{\lambda}{p\xi}D.
\end{eqnarray}
Further, from Theorem \ref{ThmMain2} we have that
\begin{eqnarray}
\P(\tau^\zeta<\infty)=\frac{\frac{\lambda}{p\xi}D}{\frac{p\xi-\lambda}{p\xi}+\frac{\lambda}{p\xi}D}.\label{expozero}
\end{eqnarray}
Representation (\ref{PRE}) given in Theorem \ref{ThmMain} and identities (\ref{exporuin}), (\ref{exponadskok}) and (\ref{expozero}) complete the proof of (\ref{expo}).
\exit

Parisian probability in this case was already identified in \cite{DassWu1}.
\\

\subsection{Brownian motion with drift}

\begin{cor}
Assume that
$$X_t=\sigma B_t+pt,$$
where $\sigma, p >0$ and $B_t$ is a standard Brownian motion.
Then
\begin{eqnarray}\label{expo2}
 \P_x(\tau^\zeta<\infty)= e^{-(2p\sigma^{-2})x} \frac{\Psi\left(\frac{p}{\sigma}\sqrt{\frac{\zeta}{2}}\right)-\frac{p}{\sigma}\sqrt{\frac{\zeta\pi}{2}}}{\Psi\left(\frac{p}{\sigma}\sqrt{\frac{\zeta}{2}}\right)+\frac{p}{\sigma}\sqrt{\frac{\zeta \pi}{2}}},
\end{eqnarray}
where $$\Psi(x)=2\sqrt{\pi}x\mathcal{N}(\sqrt{2}x)-\sqrt{\pi}x+e^{-x^2}$$
and $\mathcal{N}$(.) is a cumulative distribution function for the standard normal distribution.
\end{cor}
\proof
Note that:
$$\int_0^\infty \P(\tau_z^+>\zeta)\P_x(\tau^{-}_0<\infty, -X_{\tau^{-}_0}\in dz)=0.$$
Moreover, using
\begin{eqnarray}
\varphi(\beta)=p\beta+\frac{\sigma^2\beta^2}{2}
\end{eqnarray}
and (\ref{eq:defW})
we have
$$W(x)=\frac{1}{p}\left(1-e^{-(2p\sigma^{-2})x}\right)$$
and hence from (\ref{ruinprobref}),
\begin{eqnarray}\label{ruinabrown}
\P_x(\tau_0^-<\infty)=e^{-(2p\sigma^{-2})x}.
\end{eqnarray}
The first passage probability $\tau_\epsilon^+$ has inverse Gaussian distribution and
\begin{eqnarray}
\P(\tau_\epsilon^+<\zeta)=\mathcal{N}\left(\frac{p}{\sigma}\sqrt{\zeta}-\frac{\epsilon}{\sigma\sqrt{\zeta}}\right)+e^{2p\epsilon\sigma^{-2}}\mathcal{N}\left(-\frac{p}{\sigma}\sqrt{\zeta}-\frac{\epsilon}{\sigma\sqrt{\zeta}}\right).\label{brown1}
\end{eqnarray}
Similarly,
\begin{eqnarray}
\P_{\epsilon}(\tau_0^-<b)=e^{-2p\epsilon\sigma^{-2}}\mathcal{N}\left(\frac{p}{\sigma}\sqrt{b}-\frac{\epsilon}{\sigma\sqrt{b}}\right)+\mathcal{N}\left(-\frac{p}{\sigma}\sqrt{b}-\frac{\epsilon}{\sigma\sqrt{b}}\right).\label{brown2}
\end{eqnarray}
From (\ref{brown1})-(\ref{brown2}) using continuity of Brownian paths note that we have:
\begin{eqnarray}
\nonumber
\lim_{\epsilon\downarrow 0}\frac{p^+(b)-p(b,\zeta)}{1-p(b,\zeta)}&=&
\lim_{\epsilon\downarrow 0}\frac{\P_\epsilon(\tau_0^-<b)\P(\tau_\epsilon^+>\zeta)}{1-\P_\epsilon(\tau_0^-<b)\P(\tau_\epsilon^+\leq \zeta)}\\&=&
\frac{\sqrt{b}\Psi\left(-\frac{p}{\sigma}\sqrt{\frac{\zeta}{2}}\right)-\frac{p}{\sigma}\sqrt{\frac{b\zeta\pi}{2}}}{\sqrt{b}\Psi\left(-\frac{p}{\sigma}\sqrt{\frac{\zeta}{2}}\right)
+\sqrt{\zeta}\Psi\left(\frac{p}{\sigma}\sqrt{\frac{b}{2}}\right) }\\&=&
\frac{\sqrt{b}\Psi\left(\frac{p}{\sigma}\sqrt{\frac{\zeta}{2}}\right)-\frac{p}{\sigma}\sqrt{\frac{b\zeta\pi}{2}}}{\sqrt{b}\Psi\left(\frac{p}{\sigma}\sqrt{\frac{\zeta}{2}}\right)+\sqrt{\zeta}\Psi\left(\frac{p}{\sigma}\sqrt{\frac{b}{2}}\right) },
\end{eqnarray}
where we use fact that $\Psi(-x)=\Psi(x)$.
Note also that:
\begin{eqnarray}
\nonumber\lim_{b\to\infty}\sqrt{\frac{\zeta}{b}}\Psi\left(\frac{p}{\sigma}\sqrt{\frac{b}{2}}\right)&=&\sqrt{\frac{\zeta}{b}}\left(2\sqrt{\pi}\frac{p}{\sigma}\sqrt{\frac{b}{2}}\mathcal{N}\left(\frac{p}{\sigma}\sqrt{b}\right)-\sqrt{\pi}\frac{p}{\sigma}\sqrt{\frac{b}{2}}+e^{-(\frac{p}{\sigma})^2\frac{b}{2}}\right)
\\&=& \sqrt{\zeta}\frac{p}{\sigma}\sqrt{\frac{\pi}{2}}.
\end{eqnarray}
Putting equations (\ref{ruinabrown}) and (\ref{regular}) into (\ref{PRE}) of Theorem \ref{ThmMain} completes
the proof of (\ref{expo2}).
\exit

Probability (\ref{expo2}) was also given in \cite{DassWu1}.

One of the main goals of this paper was identifying Parisian ruin probability for more complex
spectrally negative L\'evy processes. Such typical examples we analyze below.

\subsection{General classic risk process perturbed by independent Brownian motion}
We assume that:
$$ X_t = x + p t -\sum_{i=1}^{N_t} U_i+\sigma B_t,$$
where $p, \sigma>0$ and $B_t$ is independent of Poisson process $N_t$ with intensity $\lambda$
and sequence of arriving claims $\{U_i\}_{i=1,2,\ldots}$.
Then $\varphi(\theta)=p\theta-\lambda +\lambda\int_0^\infty e^{-\theta z}\,F(dz)+\frac{\sigma^2}{2}\theta^2$
and
$\Phi(\theta)$ solves equation $$\int_0^\infty e^{-\Phi(\theta) z}\,F(dz)=\left(\lambda-p\Phi(\theta)-\frac{\sigma^2}{2}\Phi^2(\theta)+\theta\right)/\lambda.$$

Parisian ruin probability $\P_x (\tau^\zeta<\infty)$ is given in Theorem \ref{ThmMain} and all appearing there terms
except $\P(\tau^\zeta<\infty)$ are given in Remark \ref{Brown}.

Using Theorem \ref{ThmMain2}(ii) and Remark \ref{Brown2} we will find $\P(\tau^\zeta<\infty)$.
Since $X$ has a Gaussian component, then $W^{(\w)}(0+)=0$ and
from Tauberian theorem (see e.g. \cite[Theorem 1.7.1']{BGT}) it follows that
$$\lim_{\epsilon\downarrow 0}\frac{W^{(\w)}(\epsilon)}{\epsilon}=\frac{2}{\sigma^2}.$$
Thus $n(\epsilon)=\epsilon$ and
$$m(\w)=\frac{2\w}{\sigma^2\Phi(\w)}.$$
{\bf Parisian ruin probability $\P (\tau^\zeta<\infty)$ is given in
(\ref{newpzetazero}) for
\begin{equation}\label{doubleltzeta}
\int_0^\infty\int_0^\infty q(s,t)\,dt\,ds=\frac{2(\beta-\w)}{\sigma^2\beta\w(\Phi(\beta)
-\Phi \left( \w\right) )}.
\end{equation}
}

Cram\'er asymptotics $\P_x(\tau^\zeta<\infty)$ is given in Theorem \ref{Crameras} where
$\gamma$ solves equation:
$$\int_0^\infty e^{\gamma z}\,F(dz)=\left(\lambda+p\gamma-\frac{\sigma^2\gamma^2}{2}\right)/\lambda$$
and
$\widehat{\kappa}(0,0)=\varphi^\prime(0+)=p-\lambda\n$. Furthermore,
$$\mu=\lambda\int_{0}^{\infty}y e^{\gamma y} \overline{F}(y)\,dy+\sigma. $$


Assume now that claim size has
{\bf exponential distribution} $F(dz)=\xi e^{-\xi z}\,dz$.\\
Then by \cite{Wang} (see also \cite{DufGerb}):
\begin{equation}\label{ruinprobexpo}\P_x(\tau_0^-<\infty)= \widehat{\kappa}(0,0)\widehat{U}(x,\infty)=1+\sum_{i=1}^2 c_ie^{-\varrho_i x},\end{equation}
where
\begin{eqnarray*}
\varrho_1&=&\frac{(\frac{1}{2}\xi\sigma^2+p)-\sqrt{(\frac{1}{2}\xi\sigma^2+p)^2-2\sigma^2(\xi p-\lambda)}}{\sigma^2},\\
\varrho_2&=&\frac{(\frac{1}{2}\xi\sigma^2+p)+\sqrt{(\frac{1}{2}\xi\sigma^2+p)^2-2\sigma^2(\xi p-\lambda)}}{\sigma^2},\\
c_1&=&\frac{\sigma^2\varrho_2^2-2p\varrho_2}{\sigma^2(\varrho_2^2-\varrho_1^2)-2p(\varrho_2-\varrho_1)},\\
c_2&=&\frac{-\sigma^2\varrho_1^2+2p\varrho_1}{\sigma^2(\varrho_2^2-\varrho_1^2)-2p(\varrho_2-\varrho_1)}.\\
\end{eqnarray*}

We have,
\begin{eqnarray}
\lefteqn{\int_0^{\infty}e^{-\theta s}\,ds\int_0^\infty \P(\tau_z^+>s)
\P_x(\tau^{-}_0<\infty, -X_{\tau^{-}_0}\in dz)}\nonumber\\&&
=\frac{1}{\theta}\P_x(\tau_0^-<\infty)\left(1-\xi\int_0^\infty e^{-(\Phi(\theta)+\xi)s}\,ds\right)
=\frac{1}{\theta}\P_x(\tau_0^-<\infty)\left(1-\frac{\xi}{\Phi(\theta)+\xi}\right),\nonumber\\&&\label{poco}
\end{eqnarray}
where $\Phi(\theta)$ solves equation
\begin{equation}\label{phithetabr}-\frac{\sigma^2}{2}\Phi^3(\theta)-\left(\xi\frac{\sigma^2}{2}+p\right)\Phi^2(\theta)+(\lambda+\theta-p\xi)\Phi(\theta)+\xi\theta=0
\end{equation}
and by inverting this Laplace transform one can identify:
$\int_0^\infty \P(\tau_z^+>\zeta)
\P_x(\tau^{-}_0<\infty, -X_{\tau^{-}_0}\in dz).$

{\bf We will find now \begin{equation}\label{costam}
\P (\tau^\zeta<\infty)=\lim_{b\to\infty}\lim_{\epsilon\downarrow 0}\frac{(1-p(b,\zeta))-\P_\epsilon (\tau_0^->b)}{1-p(b,\zeta)}\end{equation}}
for the exponential claim size with intensity $\xi$. Note that by (\ref{poco}) Laplace transform
$\frac{1}{\theta}\left(1-\frac{\xi}{\Phi(\theta)+\xi}\right)$
is then the only single Laplace transform that must be inverted at point $s=\zeta$ to derive Parisian ruin probability $\P_x(\tau^\zeta<\infty)$
(compare with (\ref{doubleltzeta}) where double Laplace transform must be inverted).

Note that
\begin{eqnarray*}
\lefteqn{p(b,\zeta)=(1-\P_\epsilon(\tau_0^-<b, X_{\tau_0^-}=0))\xi\int_0^\infty\P(\tau_{z+\epsilon}^+\leq \zeta)e^{-\xi z}\,dz}\\&&
+\P(\tau_{\epsilon}^+\leq \zeta) \P_\epsilon(\tau_0^-<b, X_{\tau_0^-}=0)\end{eqnarray*}
and
\begin{eqnarray}
\lefteqn{\P_\epsilon(\tau_0^-<b, X_{\tau_0^-}=0)=\P_\epsilon(\tau_0^-<\infty, X_{\tau_0^-}=0)}\nonumber\\&&
-\int_0^\infty \P_z(\tau_0^-<\infty, X_{\tau_0^-}=0)
\P_\epsilon(\tau_0^-\geq b, X_{b}\in dz)\label{pass0}.\end{eqnarray}
Thus from (\ref{pass0}) and (\cite[Therorem 5.9, p. 122]{Kbook} we derive:
\begin{equation}\label{zeropassge2}
\lim_{\epsilon\downarrow 0}\frac{1-\P_\epsilon(\tau_0^-<b, X_{\tau_0^-}=0)}{\epsilon}=
d\widehat{u}^\prime(0+)-d\int_0^\infty \widehat{u}(z)
g_1(b,dz),
\end{equation}
where
\begin{equation}\label{gdz}
g_1(b,dz)=\lim_{\epsilon\downarrow 0}\frac{\P_\epsilon(\tau_0^-\geq b, X_{b}\in dz)}{\epsilon}.
\end{equation}
If
$$\lim_{\epsilon\downarrow 0}\frac{\P(\tau_\epsilon^+>\zeta)}{\epsilon}=g_2,\qquad
\lim_{\epsilon\downarrow 0}\frac{\P_\epsilon(\tau_0^->b)}{\epsilon}=g_3(b)$$
then from (\ref{costam})
$$\P (\tau^\zeta<\infty)=\frac{g_2+g_4(\infty)-g_3(\infty)}{g_2+g_4(\infty)},$$
where
$$g_4(b)=d(\widehat{u}^\prime(0+)-\int_0^\infty \widehat{u}(z)
g_1(b,dz))\xi\int_0^\infty\P(\tau_{z}^+\leq \zeta)e^{-\xi z}\,dz
$$
and
$g_i(\infty)=\lim_{b\to\infty}g_i(b)$ ($i=1,3,4$).
The probability $\P(\tau_{z}^+\leq \zeta)$ appearing in $g_4$ could be found
using Kendall's identity saying that if spectrally negative $X_t$ has has a density $m(t, x)$ at $x$, then $\tau^+_x$ also has a density at $s$ and
\begin{equation}\label{Kendall}
\frac{\P(\tau^+_x\in ds)}{ds}=\frac{x}{s}m(s,x),
\end{equation}
where in our case
$$m(s,x)=\int_{x-ps}^\infty\frac{1}{\sigma\sqrt{s}}\phi\left(s,\frac{y}{\sigma\sqrt{s}}\right)\sum_{k=1}^\infty\frac{(\lambda s)^k}{k !}e^{-\lambda s}\frac{\xi^k(y-x)^{k-1}e^{-\xi(y-x)}}{(k-1)!}.
$$
We will find functions $g_i$ for $i=1,2,3$.\\
Function $\mathbf{g_1(b,dz)}$.
Note that:
\begin{eqnarray}
\lefteqn{\P_\epsilon(\tau_0^->b, X_b\in dy)=\P(\tau_\epsilon^{B-}>b, \sigma B_b-pb\in dy)e^{-\lambda b}}\nonumber\\
&&+\lambda\int_0^b e^{-\lambda t}\int_{-\infty}^\epsilon \P(\tau_\epsilon^{B-}>t, \sigma B_t-pt\in dz) \nonumber\\&&\int_0^{-z+\epsilon} \P_{-z-h+\epsilon}(\tau^-_0>b -t, X_{b-t}\in dy)\,F(dh)\,dt,\nonumber\\&&
\label{tauzetaminus2}
\end{eqnarray}
where
$$\tau^{B-}_\epsilon=\inf\{t\geq 0: \sigma B_t-pt>\epsilon\}$$
and probability $\P_z(\tau_0^->s, X_s\in dy)$ appearing in (\ref{tauzetaminus2}) could be derived using so-called ballot theorem:
\begin{eqnarray}
\P_z(\tau_0^->s, X_s\in dy)/dy &=&m(s,y-z)-
y \int_0^s \frac{1}{s-t} m(s-t, y)m(t,-z)\, dt;\nonumber\\&&\label{main}
\end{eqnarray}
see \cite{Borovkov}. Moreover,
by the strong Markov property,
$$\P(\tau_\epsilon^{B-}>t, \sigma B_t-pt\in dz)=dz\left(\phi\left(\frac{z+pt}{\sqrt{t}\sigma}\right)-\int_0^t\P(\tau^{B-}_\epsilon\in ds)
\phi\left(\frac{z-\epsilon-p(t-s)}{\sqrt{t-s}\sigma}\right)\right),$$
where by Kendall's identity (\ref{Kendall}),
\begin{equation}\label{densityig}\P(\tau^{B-}_\epsilon\in ds)=\frac{\epsilon\, ds}{\sigma\sqrt{2\pi s^3}}\exp\left\{-\frac{(sp+\epsilon)^2}{2\sigma^2s}\right\}.
\end{equation}
Hence
\begin{equation}\label{as2}
\lim_{\epsilon\downarrow 0}\frac{\P(\tau_\epsilon^{B-}>t, \sigma B_t-pt\in dz)}{\epsilon dz}=\Xi^-(z,t),
\end{equation}
where
$$\Xi^-(z,t)=\frac{1}{2\pi}e^\frac{-p^2t+2zp}{2\sigma^2}\int_0^t\frac{1}{\sqrt{s^3}}e^{-\frac{z^2}{2\sigma^2(t-s)}}\,ds.$$
Summarizing,
\begin{eqnarray*}
\lefteqn{g_1(b,dy)=e^{-\lambda b}\Xi^-(y,b)\;dy}\nonumber\\
&&+\lambda\int_0^b e^{-\lambda t}\int_{-\infty}^0 \Xi^-(z,b)\;dz \nonumber\\&&\int_0^{-z} \P_{-z-h}(\tau^-_0>b -t, X_{b-t}\in dy)\,F(dh)\,dt.
\end{eqnarray*}

Function $\mathbf{g_3(b)}$. Similarly,
\begin{eqnarray}
\lefteqn{g_3(b)=e^{-\lambda b}\int_{-\infty}^0\Xi^-(y,b)\;dy}\nonumber\\
&&+\lambda\int_0^b e^{-\lambda t}\int_{-\infty}^0 \Xi^-(z,b)\;dz \int_0^{-z} \P_{-z-y}(\tau^-_0>b -t)\,F(dy)\,dt.\nonumber
\end{eqnarray}

Quantity $\mathbf{g_2}$. Finally,
\begin{eqnarray}
\lefteqn{\P(\tau_\epsilon^+>\zeta)=\P(\tau_\epsilon^{B+}>\zeta)e^{-\lambda \zeta}}\nonumber\\
&&+\lambda\int_0^\zeta e^{-\lambda t}\int_{-\infty}^\epsilon \P(\tau_\epsilon^{B+}>t, \sigma B_t+pt\in dz) \int_0^\infty \P_{z-y}(\tau^+_\epsilon>\zeta -t)\,F(dy)\,dt,\nonumber\\&&\label{tauzetaplus}
\end{eqnarray}
where
$$\tau^{B+}_\epsilon=\inf\{t\geq 0: \sigma B_t+pt>\epsilon\}$$
has the inverse Gaussian distribution:
$$\P(\tau^{B+}_\epsilon\in ds)=\frac{\epsilon\, ds}{\sigma\sqrt{2\pi s^3}}\exp\left\{-\frac{(sp-\epsilon)^2}{2\sigma^2s}\right\}.$$
Thus
\begin{equation}\label{as1}
\lim_{\epsilon\downarrow 0}\frac{\P(\tau_\epsilon^{B+}>t, \sigma B_t+pt\in dz)}{\epsilon dz}=\Xi^+(z,t),
\end{equation}
where
$$\Xi^+(z,t)=\frac{1}{2\pi}e^\frac{-p^2t-2zp}{2\sigma^2}\int_0^t\frac{1}{\sqrt{s^3}}e^{-\frac{z^2}{2\sigma^2(t-s)}}\,ds$$
and
\begin{eqnarray}
\lefteqn{g_2=e^{-\lambda \zeta}\int_{-\infty}^0\Xi^+(y,\zeta)\;dy}\nonumber\\
&&+\lambda\int_0^\zeta e^{-\lambda t}\int_{-\infty}^0 \Xi^+(z,\zeta) \int_0^\infty \P_{z-y}(\tau^+_0>\zeta -t)\,F(dy)\,dt.
\end{eqnarray}

The probability $\P_{z}(\tau^+_0>s)=\P(\tau_{-z}^+>s)$ for $z<0$ appearing above could be obtained using Kendall's identity (\ref{Kendall}).

\subsection{General classic risk process perturbed by $\alpha$-stable motion}
We will assume now that:
\begin{equation}\label{rstablep}
 X_t = x + p t -\sum_{i=1}^{N_t} U_i+ Z_t,\end{equation}
where $p, \sigma>0$ and $Z_t$ is a spectrally negative
$\alpha$-stable motion with $\alpha\in (1,2)$
independent of Poisson process $N_t$ with intensity $\lambda$
and the sequence of arriving claims $\{U_i\}_{i=1,2,\ldots}$.
Note that for $\alpha \in (1,2)$ process $Z$ is of unbounded variation. Moreover, process $X$ has no gaussian component and it does not creeps downward.
Then $\varphi(\theta)=p\theta-\lambda +\lambda\int_0^\infty e^{-\theta z}\,F(dz)+c\theta^\alpha$
for some $c>0$ and
$\Phi(\theta)$ solves equation $$\int_0^\infty e^{-\Phi(\theta) z}\,F(dz)=\left(\lambda-p\Phi(\theta)-c\Phi^\alpha(\theta)+\theta\right)/\lambda.$$
Then from \cite{furrer} we have:
\begin{equation}\label{stable}
\P_x(\tau_0^-<\infty)=\widehat{\kappa}(0,0)\widehat{U}(x,\infty)
=1-(1-\rho_0)\sum_{n=0}^\infty\rho_0^n\left(K^{(n+1)*}*M^{n*}\right)(x),
\end{equation}
where $\rho_0=\frac{\lambda\nu}{p}$, $M(dx)=\frac{1}{\nu}\overline{F}(y)\,dy$ and $\overline{K}(dx)=
\sum_{n=0}^\infty \frac{(-p)^n}{\Gamma(1+(\alpha-1)n}x^{(\alpha-1)n}$.

Finally, {\bf $\P_x(\tau^\zeta<\infty)$ could be found} using main representation given Theorem
\ref{ThmMain}, the equation (\ref{ltx0}) and identity (\ref{stable}).
Using Theorem \ref{ThmMain2}(ii) and Remark \ref{Brown2} {\bf one can also identify
$\P(\tau^\zeta<\infty)$} where
similarly like in case of Brownian perturbation we have:
$W^{(\w)}(0+)=0$ and
$$\lim_{\epsilon\downarrow 0}\frac{W^{(\w)}(\epsilon)}{\epsilon^{\alpha-1}}=\frac{1}{c\Gamma(\alpha)}.$$
Thus $n(\epsilon)=\epsilon^{\alpha -1}$ and
$$m(\w)=\frac{\w}{c\Gamma(\alpha)\Phi(\w)}.$$

Note also that for process (\ref{rstablep}) we have
\begin{equation}\label{stablejump}\Pi_{\widehat{X}}(dy)=\frac{c}{y^{1+\alpha}}\,dy+\lambda F(dy).\end{equation}
Assuming that $\overline{F}(x)={\rm o}(x^{-\alpha})$ for large $x$ we derive that
$\overline{\Pi}_{\widehat{X}}\in\mathcal{S}^{(0)}$.
From Theorem \ref{Thmconveq} we have
$$\P_x(\tau^\zeta<\infty)\sim \left(\P(\tau^\zeta<\infty)/\E X_1 +\P(\tau^\zeta=\infty)f^{(e)}(\zeta)\right)\frac{c}{\alpha(\alpha-1)}x^{-\alpha+1},
$$
as $x\to\infty$, where the Laplace transform of $f^{(e)}(\cdot)$ is given in (\ref{fe}) for $ \Pi_{\widehat{X}}$
defined in (\ref{stablejump}).

\section*{Acknowledgements}
This work is partially supported by the Ministry of Science and
Higher Education of Poland under the grants N N201 394137
(2009-2011) and N N201 525638 (2010-2011).

\end{document}